\documentstyle[11pt]{article}
\catcode`\@=11
\def\[#1{\relax\ifmmode\@badmath\else\begin{equation}\label{#1}\fi}
\def\]{\relax\ifmmode\ifinner\@badmath\else\end{equation}\fi
         \else \@badmath \fi}
\newcounter{theorem}[section]
\def\thetheorem{\thesection.\arabic{theorem}}

\def\proclaim{\pr@claim{\bf}{ \thetheorem}{\sl}}
\def\heading{\pr@claim{\bf}{ \thetheorem}{\relax}}

\outer\def\conjecture #1.{\proclaim Conjecture #1.}
\outer\def\theorem #1.{\proclaim Theorem #1.}
\outer\def\proposition #1.{\proclaim Proposition #1.}
\outer\def\lemma #1.{\proclaim Lemma #1.}
\outer\def\corollary #1.{\proclaim Corollary #1.}
\outer\def\remark #1.{\pr@claim{\sl}{ \thetheorem}{\relax}Remark #1.}
\outer\def\definition #1.{\heading Definition #1.}
\outer\def\acknowledgement{\pr@claim{\sl}{\relax}{\relax}Acknowledgement *.}

\def\pr@claim#1#2#3#4 #5.#6\par{\refstepcounter{theorem}\medbreak\noindent
       \def\next{#5}{#1#4#2\if\next*\else\label{\next}\fi.}\nobreak
       {\enspace#3#6}\par
  \ifdim\lastskip<\medskipamount \removelastskip\penalty55\medskip\fi
  }

\def\qedsymbol{$\Box$}

\def\Bbb{\ifmmode\let\next\Bbb@\else
 \def\next{\errmessage{Use \string\Bbb\space only in math mode}}\fi\next}
 \def\Bbb@#1{{\Bbb@@{#1}}}
 \def\Bbb@@#1{\rm\bf#1}              

\def\itm#1{\par\hang\textindent{#1}}

\newif\ifcomment
\def\comment{\ifcomment\bgroup\par\medskip\noindent\small}
\def\endcomment{\par\medskip\noindent\egroup}

\newif\ifnote

\let\\=\cr

\def\mathalign#1#2#3#4{\null\,\vcenter\bgroup\openup\jot\m@th
\ialign\bgroup\strut#1$\displaystyle##$#2&&#3$\displaystyle{}##$#4\crcr}
\def\endmathalign{\crcr\egroup\egroup\,}

\def\txt#1{\hbox{ #1 }}

\outer\def\proof#1.{\smallbreak\noindent{\sl Proof#1\/}.\ }

\def\endproof{\removelastskip\nobreak
{\unskip\nobreak\hfil\penalty50            
\hskip3em\hbox{}\nobreak\hfil\qedsymbol
\parfillskip=0pt \finalhyphendemerits=0 \par
\medbreak\noindent\ignorespaces}}

\newdimen\arrowlen
\def\m@p#1#2#3#4{
       \buildrel
              \hbox spread \arrowlen{\skip@=-0.33\arrowlen plus 1 fil
                    \hskip\skip@$\m@th\scriptstyle #4$\hskip\skip@
              }
       \over{
              \mathord#1\mkern-6mu
              \cleaders\hbox{$\mkern-2mu\mathord#2\mkern-2mu$}\hfill
              \mkern-6mu\mathord#3
       }
}

\def\map{\m@p--\rightarrow}

\def\rmmath#1{\mathop{\rm #1}\nolimits}

\def\NS{\rmmath{NS}}

\let\slant=/

\mathcode`\:="603A 
\mathchardef\:="303A 

\def\C{{\Bbb C}}
\def\Q{{\Bbb Q}}
\def\Z{{\Bbb Z}}

\def\O{{\cal O}}

\def\numfrac#1#2{\mathchoice{{\textstyle{ #1\over#2}}}%
{{ #1\over#2}}{{#1/#2}}{{#1/#2}}}

\def\half{{\numfrac12}}

\let\next=\~
\let\~=\tilde
\let\tilde=\next

\let\next=\^
\let\^=\hat
\let\hat=\next

\def\"#1{{\accent"7F \if#1i\i\else#1\fi}}
\def\({\left(}
\def\){\right)}

\def\SU{{\rm SU}}

\def\M{{\cal M}}
\def\<#1>{\left<#1\right>}

\catcode`\@=12

\title{Some remarks on the Kronheimer-Mrowka classes of algebraic surfaces}
\author{Rogier Brussee%
\thanks{Mathematisches Institut Universit\"at Bayreuth,
Universit\"ats Stra\ss e 30, D-98440 Bayreuth, Germany,
e-mail:Brussee@btm8x1.mat.uni-bayreuth.de}%
}
\date{alg-geom/9308003} 
\def\KM{KM-}
\def\NE{\rmmath{NE}}
\def\NEbar{\rmmath{\overline{NE}}}
\def\NS{\rmmath{NS}}
\def\ddz#1#2{{\partial^{#1} #2  \over \partial z^{#1}}}
\def\SO{{\rm SO}}

\begin{document}

\maketitle

\section{Introduction}

Recently Kronheimer and Mrowka have announced a very interesting result,
which sheds new light on the Donaldson polynomials \cite{\KM}.
They find recurrence relations between the Donaldson polynomials,
by finding relations between the polynomials and the minimal genus
of a smooth real surface representing an homology class.
To be more precise we need a definition.

For a simply connected 4-manifold $X$ with odd $b_+ \ge 3$, we denote
the $\SU(2)$ polynomials on $H_0(X) \oplus H_2(X)$ by $q_k(X)$. $X$ is
called {\sl simple} if we have
$$
      q_k(X)(pt^2,-) = 4 q_{k-1}(X),
                \qquad d= 4k -\numfrac32(1+b_+).
$$
For simple 4-manifolds it is convenient to label the polynomials by their
degree on $H_2(X)$ i.e. we define
$$
       q_d(X) =\cases{q_k(X)|_{H_2(X)} &if $d = 4k  -\numfrac32(1 + b_+)$,
\\
             q_k(X)(pt,-)|_{H_2(X)}&if $d = 4k - 2 -\numfrac32(1 +b_+)$,
\\
             0                      & otherwise}
$$
The
{\em Donaldson series} is then the formal power series
$q(X) = \sum_d q_d(X) / d!\,$.

\theorem KM. (Kronheimer, Mrowka)
For every simple 4-manifold $X$
there exist a finite number of {Kronheimer-Mrowka classes}
$K_1,\ldots, K_p\in H^2(X)$ and
non zero rational numbers $a_1,\ldots a_p$ such that
\itm{(i)} $q(X) = e^{Q / 2} \sum_{i = 1}^n  a_i e^{K_i}$,
\itm{(ii)} $K_i \equiv w_2(X) \pmod 2$, for all $i=1,\ldots,p$,
\itm{(iii)} if $K_i \in \{K_1,\ldots, K_p\}$ then
           $-K_i\in \{K_1,\ldots K_p\}$,
\itm{(iv)} for every homologically nontrivial connected real surface
    $\Sigma$ with $\Sigma^2 \ge 0$ and every Kronheimer-Mrowka
    class $K_i$ we have
$$
       2g(\Sigma) - 2 \ge \Sigma^2 + K_i\cdot \Sigma.
$$

Here $Q$ is the intersection form. The Kronheimer Mrowka classes will
be abbreviated KM-classes and
are called the
basic classes in \cite{\KM}. Condition (ii) is    reminiscent
of the Wu formula, whereas condition (iv) is similar to the
the adjunction  formula for the
genus of a smooth algebraic curve. This suggests that  for complex surfaces,
the  \KM classes should be closely related to the canonical divisor
$K_X$. Indeed in the examples  and the conjectural expression for
the Donaldson series of elliptic surfaces in the announcement,
the \KM classes are of type $K_i = \alpha_i K_{\min} + \sum
\beta_{ij} E_j$ where $K_{\min}$ is the canonical divisor of the
minimal model, $E_1,\ldots E_l$ the $(-1)$-curves, $\alpha_i$ a
rational number with $|\alpha_i| \le 1$ and $\beta_{ij} = \pm 1$. A
formulation which is a little less obvious but which
generalises better, as we will see, is $K_i = C-D$
where $C$, $D$  are  divisors such that $K_X = C+D$ and such that a
multiple is effective. Hence in these cases the canonical class is a \KM class
which is extremal from an algebraic geometric point of view.


To formulate the general relationship, we recall that
the effective cone  of a complex surface $\NE(X)$,
is the positive rational cone in
$\NS(X)_\Q \subset H^2(X,\Q)$ generated by effective divisors.
Let $\NEbar(X)$ be its closure in  the norm topology.

\theorem main.
For every \KM class $K_i$ on a simple simply
connected algebraic surface $X$, there is a unique decomposition
$K_X = C_i + D_i$ in  $\Z$-divisors $C_i,D_i \in \NEbar(X)$ such that
$ K_i = C_i - D_i $. In particular $K_i = K_X$ if and only if there is
a smooth hyperplane section $H$ such that $2g(H) -2 = H^2 +K_i\cdot H$.

The results in \cite{Kronheimer:genus} seem to indicate that a \KM class
and a hyperplane section $H$ as in the theorem exist, at least
when there is an
$\omega \in H^0(K)$ such that for sufficiently large $k$,
$q_k(\omega + \bar\omega) \ne 0$. For minimal surfaces of general type this
would imply the invariance of the canonical class up to sign under orientation
preserving self-diffeomorphisms.

\corollary gentype.
Assume in addition that $X$ is minimal and of general type, then $K_i^2 \le
K_X^2$ with equality if and only if $K_i =\pm K_X$.


By the Lefschetz (1,1) theorem \cite[p. 163]{G&H},
the algebraicity of the \KM classes is equivalent to the $K_i$ being of type
$(1,1)$,
In fact this is what we will prove, using
that the Donaldson polynomials are of pure Hodge type
as in \cite{(-1)-curve}.
Since we assume that $p_g >0$ this shows that the lattice
spanned by the $K_i$ is  a proper sublattice.
Moreover since the $K_i$ are defined by the differentiable structure,
they are contained in the fixed lattice of the variation of
Hodge structures defined by a family of complex
structures on $X$. In favourable circumstances  this should force
the \KM classes  to be
in $H^2(X,\Z) \cap [K_{\min}, E_1,\ldots, E_l]$.

Another implication of the theorem is that the
\KM classes are trivial on  $H^{0,2}(X)  \subset H^2(X,\C)$. Thus we get
the following corollary.

\corollary forms.
If $X$ is a simple and simply connected surface, then
for all $\omega \in H^0(K_X)$ we have
$$
       q(\omega + \bar\omega) = q_0 e^{\int \omega \wedge \bar\omega}
$$
where $q_0$ is the Donaldson polynomial of degree 0

It would be rather interesting to understand this formula from an
algebraic geometric point of view, possibly clarifying the role of the
simple\-ness condition. Combining the corollary with O'Grady's non
vanishing
result, we see that $q_0 \ne 0$ if $X$ is  of general type, $p_g$ is odd
and $|K_{\min}|$, the linear system of the canonical class of the minimal
model, contains a reduced curve \cite[th. 2.4]{O'Grady},
\cite[appendix]{JunLi:Kodaira}, \cite[th. 1]{Morgan}.

In a similar vein, since the Neron-Severi group has a non degenerate
intersection form, the $K_i$ are determined by their intersection
products with divisors. Thus we get

\corollary alg.
For a simple simply connected surface, the Donaldson series
$q(X)$ is determined by  $q(X)|_{\NS(X)}$.

The corollary says that by knowing the algebraic part of all Donaldson
polynomials we can reconstruct the transcendental part as well, i.e. in
the simple case, the polynomials defined by Jun-Li
\cite{JunLi:polynomials} contain as much information as the full
polynomials. Moreover since $\NEbar \cap K_X -\NEbar$ is a
bounded subset of $\NS(X,\Q)$, and the $K_i$ are integral,
at least in principle, we get
an effective bound on the number of Kronheimer Mrowka classes,
hence on the number of
polynomials one has to compute in order to reconstruct the
Donaldson series.

\section{Proof of theorem \protect\ref{main} and
corollary \protect\ref{gentype}}

We have to prove that the \KM classes are of  type $(1,1)$. Accepting this,
 the rest of the statement of theorem \ref{main} and corollary \ref{gentype}
 is a consequence of property (ii), (iii), and (iv) of the \KM
classes.

By property (ii) we can write $K_i = K_X - 2D$ for some $\Z$-divisor $D$.
For every very ample line bundle $\O(H)$, we choose a smooth
connected hyperplane section $H$. Then we get
$$
       2g(H)-2 = H^2 + H\cdot K_X \ge H^2 + H\cdot K_i
                                  = H^2 +H\cdot K_X - 2H \cdot D
$$
i.e. $D \cdot H \ge 0$. Thus by the duality of the closure of
the effective cone $\NEbar(X)$ and the nef cone, we conclude that
$D \in \NEbar(X)$ \cite[prop. 2.3]{Wilson:birational}.
We also have $ C:= (K_X + K_i)/2 \in \NEbar$ by property (iii).
Rewriting we get $K_X = C+D$ and $K_i= C-D$ as claimed.
Note that nothing is gained by applying the inequality to other
smooth connected divisors $C$ with $C^2 \ge 0$, since such divisors are nef.
Finally note that $H\cdot K_i = 2g(H)-2-H^2 = H\cdot K_X$ if and only if
$D = 0$, since $\NEbar \cap H^\perp = (0)$.
This proves theorem \ref{main} up to algebraicity.

Now assume temporarily that $X$ is minimal  and of general type.
Write $K_i^2 = K_X^2 + 4(D^2 -K_X \cdot D)$. Since $K_X$ is nef, we get
$K_i^2 \le K_X^2$ if $D^2 \le 0$, with equality iff $D^2 = K\cdot D = 0$.
The latter is equivalent to $D=0$ by the Hodge index theorem.
Interchanging $K_i$ and $-K_i$ if necessary, we are left with the case
$C^2>0$, $D^2>0$. Since $C,D \in \NEbar$, we have $C\cdot D >0$ by the
numerical connectedness of $K_X$ \cite[prop VII.6.1]{BPV} (strictly
speaking we need that $C$ and $D$ are effective, but the proof of [loc. cit]
carries over without change). Thus $K_X^2 > D^2$. Then again by the Hodge index
theorem we get $(K\cdot D)^2 \ge K^2D^2 > (D^2)^2$, so  $K_i^2 < K^2$.
This proves the corollary.

To prove that the \KM classes are of type (1,1),
we need a slight generalization of \cite[prop. 3.1]{(-1)-curve}. For
simplicity we restrict ourselves to the $\SU(2)$ case and a statement
about Hodge types, but the
proof can easily be modified to show that all $\SO(3)$
polynomials $q_{L,k}(X)$ with  $L \in \NS(X)$ come form an algebraic
cycle (cf. [loc. cit.]).
Consider the Donaldson polynomial $q_d$ as an element of $S^d H^2(X)$
en\-dowed with its natural Hodge structure.

\lemma pure.
For every $d \ge 0$, the Donaldson
polynomials $q_d$ are pure of type $(d,d)$.

\proof.
We temporarily index the polynomials on $H_2(X)$ by $k = c_2$.
By \cite[prop. 3.1]{(-1)-curve} the lemma is true for $q_k$ with odd
$k \gg 0$ and evaluated on 2 dimensional classes.
Here sufficiently large means that the moduli space $\M_k(X)$ of
stable bundles on $X$ for generic polarisation $H$ is generically smooth
and reduced of the proper dimension, and that the lower moduli spaces
$\M_{k'}(X)$ for $k'< k$ have sufficiently low dimension,
so that we can apply Morgan's
comparison result \cite{Morgan}. For the general case we use stabilization.

By \cite[th. 2.1.1]{Morgan-O'Grady}, for every $k$ we can find an $l_0$,
and $\epsilon_1,\ldots, \epsilon_{l}$ sufficiently small
such that for a generic polarisation on the $l \ge l_0$-fold blow-up of type
$H - \sum \epsilon_i E_i$, the moduli space
$\M_{k+l}(\^X(x_1,\ldots, x_l))$ is generically smooth of the proper
dimension and the lower moduli spaces $\M_{k'}(\^X)$ with  $k' <k$,
have sufficiently low dimension.
Thus  $q_{k+l}(\^X(x_1,\ldots x_l))$ is pure of type
$(d+4l,d+4l)$

Now by the blow-up  formula \cite[th. 4.3.1]{F&M}, we have
$$
     q_k(X) = (-\half)^l q_{k+l}(\^X(x_1,\ldots x_l))(E_1^4,\ldots,E_l^4,-).
$$
Clearly the exceptional divisors $E_i$ are pure of type $(1,1)$, hence
$q_k(X)$ is pure of type $(d,d)$. Finally
$$
       q_k(X)(pt,-) = -\half q_{k+1}(\^X(x))(E^6,-),
$$
by \cite[below cor. 4.3.2]{F&M} (without proof) or \cite{consum},
and so we can use the same argument as above.
\endproof

The rest of the proof is now straightforward. Let
$$
       C(X) = e^{-Q/2}q(X) = \sum a_i e^{K_i}.
$$
Then $C(X) = \sum C_d(X)$, where $C_d(X)$ is a homogeneous
polynomial of pure type $(d,d)$, since both $\exp(-Q/2)$ and $q(X)$ is a
sum of homogeneous polynomials of pure type.
Write $K_i = \alpha_i + \beta_i + \bar\alpha_i$, with
$\alpha_i \in H^{0,2}$, $\beta_i \in H^{1,1}$. Let $z \in H^{1,1}$ be
variable and denote by $\ddz{}{}$ the directional derivative.  Then we have
$$
       \ddz n{C(X)} = \sum_{d=0}^\infty \ddz n{C_d}
                   = \sum a_i \<\beta_i,z>^n e^{K_i}.
$$
Clearly $\ddz n{C_d}$ has pure type $(d-n,d-n)$. Thus,
upon restriction to $H^{0,2}(X)$ only the constant term contributes,
and we get
$$
       \left.\ddz n{C(X)}\right|_{H^{0,2}} = \txt{constant}
                          = \sum a_i \<\beta_i,z>^n e^{\bar\alpha_i}.
$$
We conclude that
$$
       \sum_{\alpha_i \ne 0} a_i \<\beta_i,z>^n e^{\bar\alpha_i} =0.
$$
Since $z$ and $n$ are arbitrary we find that the sum must be
empty, i.e. all \KM classes $K_i$ are of type $(1,1)$.
\endproof


\begin{thebibliography}{BPV}\frenchspacing

\def\Springer{Heidelberg, Berlin, New York: Springer }
\parskip 0pt plus \baselineskip

\bibitem[BPV]{BPV}
Barth, W., Peters, C. and Van de Ven, A.:
{\em Compact complex surfaces}.
\Springer 1984.

\bibitem[Br1]{(-1)-curve}
Brussee, R.:
{\em  On the $(-1)$ curve conjecture of Friedman and Morgan},
Inv. Math.  to appear 1993, e-print alg-geom/9209001.

\bibitem[Br2]{consum}
Brussee, R.:
{\em Donaldson polynomials of connected sums},
in preparation.

\bibitem[F-M]{F&M}
Friedman, R. and Morgan, J.W.:
{\em Smooth four manifolds and complex surfaces}.
Ergebnisse der Mathematik und ihrer Grenzgebiete 3. Folge
\Springer (to appear 1993).

\bibitem[G-H]{G&H}
Griffiths, P., and   Harris, J.:
{\em Principles of algebraic geometry}.
New York Chichester Brisbane Toronto: John Wiley \& Sons 1978.

\bibitem[Kro]{Kronheimer:genus}
Kroneimer, P.B.:
{\em The genus minimizing property of algebraic curves}.
Bull. A.M.S. {\bf 29} p. 63--69 July 1993.

\bibitem[K-M]{\KM}
Kronheimer P.B., and Mrowka T.S.:
{\em Recurrence relations and asymptotics for four-manifold invariants}.
research announcement submitted to Bull. Amer. Math. Soc. 1993.

\bibitem[Li1]{JunLi:polynomials}
Li, J.:
{\em Algebraic geometric interpretation of Donaldson's polynomial invariants}.
J. Diff. Geom. {\bf 37} 417--466 (1993).

\bibitem[Li2]{JunLi:Kodaira}
Li, J.:
{\em Kodaira dimension of moduli spaces of vector bundles of surfaces},
preprint UCLA 1993.

\bibitem[Mor]{Morgan}
Morgan, J.W.:
{\em Comparison of the Donaldson polynomial invariants with their
algebro-geometric analogues}.
preprint Columbia university New York (1992).

\bibitem[M-O]{Morgan-O'Grady}
Morgan, J., and O'Grady, K.:
{\em The smooth classification of fake K3's and similar surfaces}.
preprint 1992.

\bibitem[O'G]{O'Grady}
O'Grady, K.:
{\em Algebro-geometric analogues of Donaldson's polynomials}.
Inv. Math. {\bf 107} 351--359 (1992).

\bibitem[Wil]{Wilson:birational}
Wilson, P.:
{\em Towards  birational classification of algebraic varieties}
Bull. London Math. Soc. {\bf 19}  1--48 (1987)
\end{thebibliography}
\end{document}